\documentclass[12pt,a4paper]{amsart}
\usepackage[all]{xy}

\usepackage{diagrams,amssymb}
\diagramstyle[height=2em,width=2em,midshaft,labelstyle=\scriptstyle]
\newarrow{To}----{->}
\newarrow{Epi}----{->>}
\newarrow{Mono}>---{->}
\newarrow{Iso}>---{->>}
\newarrow{Mapsto}|---{->}
\newarrow{Igual}=====
\newarrow{Dashto}{}{dash}{}{dash}{->}

\newcommand{\Z}{{\mathbb Z}}

\newcommand{\F}{{\mathbb F}}

\newcommand{\pcom}{{}_{p}^{\wedge}}
\newcommand{\doscom}{_{2}^{\wedge}}

\DeclareMathAlphabet\EuR{U}{eur}{m}{n}
\SetMathAlphabet\EuR{bold}{U}{eur}{b}{n}

\newcommand{\Ext}{\operatorname{Ext}\nolimits}

\newcommand{\defeq}{\overset{\text{\textup{def}}}{=}}

\newcommand{\cala}{\mathcal{A}}

\newcommand{\calf}{\mathcal{F}}

\newcommand{\call}{\mathcal{L}}

\newcommand{\widebar}[1]{\overset{\mskip3mu\hrulefill\mskip3mu}{#1}
                \vphantom{#1}}

\let\oldcirc=\circ
\renewcommand{\circ}{\mathchoice
    {\mathbin{\scriptstyle\oldcirc}}{\mathbin{\scriptstyle\oldcirc}}
    {\mathbin{\scriptscriptstyle\oldcirc}}
    {\mathbin{\scriptscriptstyle\oldcirc}}}

\newcommand{\hclim}[1]{\setbox1=\hbox{\rm hocolim}
    \setbox2=\hbox to \wd1{\rightarrowfill} \ht2=0pt \dp2=-1pt
    \mathop{\vtop{\baselineskip=5pt\box1\box2}}
    _{#1}}

\newcommand{\higherlim}[2]{\displaystyle\setbox1=\hbox{\rm lim}
	\setbox2=\hbox to \wd1{\leftarrowfill} \ht2=0pt \dp2=-1pt
	\setbox3=\hbox{$\scriptstyle{#1}$}
	\ifdim\wd1<\wd3
	\mathop{\hphantom{^{#2}}\vtop{\baselineskip=5pt\box1\box2}^{#2}}_{#1}
	\else
	\mathop{\vtop{\baselineskip=5pt\box1\box2}}\limits_{#1}\nolimits^{#2}
	\fi}

                     \newcommand{\Hom}{\operatorname{Hom}\nolimits}

\newcommand{\Ker}{\operatorname{Ker}\nolimits}

\renewcommand{\Im}{\operatorname{Im}\nolimits}

\newcommand{\xxto}[1]{\mathrel{\mathop{%
  \setbox0\hbox{$\ {\scriptstyle#1}\ $}%
  \hbox to \wd0{\rightarrowfill}}^{#1}}%
}

\newcommand{\xto}[2][]{%
  \mathrel{\mathop{%
    \setbox0\vbox{
      \hbox{$\scriptstyle\;\;{#1}\;\;$}%
      \hbox{$\scriptstyle\;\;{#2}\;\;$}%
    }%
    \hbox to\wd0{\rightarrowfill}\displaystyle}%
  \limits^{#2}\ifx{#1}{}\else{_{#1}}\fi}%
}

\newcommand{\longleft}[1]{\;{\leftarrow%
\count255=0 \loop \mathrel{\mkern-6mu}%
    \relbar\advance\count255 by1\ifnum\count255<#1\repeat}\;}
\newcommand{\longright}[1]{\;{\count255=0 \loop \relbar\mathrel{\mkern-6mu}%
    \advance\count255 by1\ifnum\count255<#1\repeat\rightarrow}\;}

\newcommand{\RIGHT}[3]{\mathrel{\mathop{\kern0pt\longright#1}
        \limits^{#2}_{#3}}}

\newcommand{\LEFT}[3]{\mathrel{\mathop{\kern0pt\longleft#1}\limits^{#2}_{#3}}
}
\newcommand{\dRIGHT}[3]{\mathrel{%
   \mathop{\vcenter{\baselineskip=0pt\hbox{$\kern0pt\longright#1$}%
   \hbox{$\kern0pt\longright#1$}}}\limits^{#2}_{#3}}}
\newcommand{\LRIGHT}[3]{\mathrel{%
   \mathop{\vcenter{\baselineskip=0pt\hbox{$\kern0pt\longleft#1$}%
   \hbox{$\kern0pt\longright#1$}}}\limits^{#2}_{#3}}}
\newcommand{\RLEFT}[3]{\mathrel{%
   \mathop{\vcenter{\baselineskip=0pt\hbox{$\kern0pt\longright#1$}%
   \hbox{$\kern0pt\longleft#1$}}}\limits^{#2}_{#3}}}
\newcommand{\onto}[1]{\;{\count255=0 \loop \relbar\joinrel
    \advance\count255 by1
    \ifnum\count255<#1 \repeat \twoheadrightarrow}\;}

\newtheoremstyle{slant}{}{}{\slshape}{}{\bfseries}{.}{.5em}{}%
\newtheoremstyle{special}{}{}{\slshape}{}{\bfseries}{.}{.5em}{\thmnote{#3}}

\newtheorem{Thm}{Theorem}[section]
\newtheorem{Prop}[Thm]{Proposition}

\newtheorem{Th}{Theorem}

\theoremstyle{remark}

\newcommand{\Cotor}{\mathrm{Cotor}}

\title{Loop space homology associated to the mod 2 Dickson invariants}
\author{Ran Levi}
\address{Institute of Mathematics, University of Aberdeen,
Aberdeen AB24 3UE, U.K.}
\email{r.levi@abdn.ac.uk}
\author{Nora Seeliger}
\address{L.A.G.A. -- Laboratoire de Math\'ematiques, Institut Galil\'ee, 99 Avenue. J.B. Cl\'ement, 93430
Villetaneuse, France}
\email{seeliger@math.univ-paris13.fr}


\newcommand{\Sol}{\mathrm{Sol}}

\begin{document}
\maketitle

\begin{abstract}
The spaces $BG_2$ and $BDI(4)$ have the property that their mod 2 cohomology is given by the rank 3 and 4 Dickson invariants respectively. Associated with these spaces one has for $q$ odd the classifying spaces of the finite groups $BG_2(q)$ and the exotic family of classifying spaces of 2-local finite groups $B\Sol(q)$. In this article we compute the loop space homology of $BG_2(q)\doscom$ and $B\Sol(q)$ for all odd primes $q$, as algebras over the Steenrod algebra, and the associated Bockstein spectral sequences. 
\end{abstract}

It is well known that the mod 2 Dickson invariants $P[u_1,u_2,\ldots, u_n]^{GL_n(\F_2)}$ are realisable as the mod 2 cohomology of a space only for $n\le 4$. For $n=2,3$ the corresponding spaces are the classifying spaces of the Lie groups  $SO(3)$ and $G_2$ respectively. For $n=4$ a space $BDI(4)$ was constructed by Dwyer and Wilkerson, which realises the rank 4 invariants \cite{DW}. In 1994  Benson introduced a family of spaces $B\Sol(q)$, one for each odd prime power $q$, closely related to $BDI(4)$, which he claimed realised the exotic fusion patterns studied by Solomon 20 years earlier \cite{Benson, Sol}. He obtained this family of spaces by considering the pullback of the system
\[BDI(4)\rTo{\psi^q\top1} BDI(4)\times BDI(4)\lTo{\Delta} BDI(4),\]
where $\psi^q$ is the degree $q$ unstable Adams operation constructed by Notbohm \cite{Notbohm}. In \cite{LO} the first named author and Oliver showed that the patterns studied by Solomon form in fact what became known more recently as saturated fusion systems, and that these fusion systems admit associated centric linking systems, and thus give rise to a family of 2-local finite groups (see \cite{BLO2}). The "classifying spaces" of these 2-local finite groups are also named $B\Sol(q)$, and are shown to coincide with Benson's family. The family $B\Sol(q)$ provides one of the most interesting collections of p-local finite groups, in that they are all exotic, and  to date the only exotic systems known at the prime 2. The module structure of $H^*(B\Sol(q),\F_2)$ was calculated by Benson, and the algebra and $\cala_2$-module structure were determined by Grbic, who also computed the Bockstein spectral sequence for these spaces  \cite{Grbic}.

In this article we consider the spaces $BG_2(q)\doscom$ and $B\Sol(q)$ for all odd prime powers $q$, and present a complete  calculation of their loop space homology. There are strong results known on the homotopy type of $\Omega BG\pcom$ when $G$ is a finite group \cite{Levi}, but not much is known on loop spaces of exotic classifying spaces. Furthermore, as we shall see these two families exhibit very systematic behaviour, which might be worth exploring further. This motivates our calculations.

Throughout this paper $H_*(-)$ and $H^*(-)$ will mean mod 2 homology and cohomology respectively. Different
coefficients will always be explicitly specified. Subscripts on homology or cohomology classes will always denote their degrees. The letters $P$, $E$, $T$ and $\Gamma$ will be used to denote the polynomial, exterior, tensor and divided power algebras respectively. By convention we will always use the notation $T[x]$ to denote a tensor algebra on a single odd-dimensional generator, although over $\F_2$ the  tensor algebra on a single generator in any dimension is graded commutative. 
The spectral
sequences of Serre, Bockstein and  Eilenberg-Moore  will be used in
our calculations and will be abbreviated as SSS, BSS and EMSS
respectively. 

Our first result is a calculation of the mod 2 loop space homology of $BG_2(q)$, for any odd prime power $q$. 

\begin{Th} \label{G2}
Fix an odd prime power q. Then
\[H_*(\Omega BG_2(q)\doscom)\cong P[a_2]/(a_2^2)\otimes P[a_4, b_{10}]\otimes E[x_3,x_5]\otimes P[z_6]/(z_6^2),\]
as modules over $H_*(\Omega BG_2)\cong P[a_2]/(a_2^2)\otimes P[a_4, b_{10}]$. Furthermore:
\begin{itemize}
\item The relations which determine the algebra extension are given by $x_3^2 = x_5^2 = z_6^2 =0$, $[a_2,z_6] = a_4^2$, $[a_4,z_6] = b_{10} + a_2a_4^2$, and $[b_{10}, z_6] = a_4^4$. All other commutators of generators are trivial. 
\item The reduced coproduct is given by $\widebar{\Delta}(a_4) = a_2\otimes a_2$, $\widebar{\Delta}(z_6) = x_3\otimes x_3$, while all other generators are primitive. 
\item The action of the dual Steenrod algebra is determined by 
\begin{center}
 \begin{tabular}{||l|c|c|l|l|c|l||} \hline
&  $a_2$ & $x_3$  & $a_4$  & $x_5$&  $z_6$ & $b_{10}$   \\ \hline
$Sq^1_*$  &  $0$ &$0$  &$0$ & $0$  &  $x_5$ & $0$   \\ \hline
$Sq^2_*$ &  $0$ & $0$ & $a_2$ &$x_3$ &  $0$ & $a_4^2$  \\ \hline
\end{tabular}
\end{center}
and the Steenrod axioms.
\item The homology Bockstein spectral sequences are determined by 
\begin{center}
 \begin{tabular}{||l|c|c|c|c|c|c|c|l||} \hline
$q\equiv 1(4)$& $a_2$ & $x_3$ & $a_4$ & $x_5$  & $a_2x_3$ & $z_6$ & $b_{10}$ & $x_5z_6$\\ \hline
$Sq^1_*$  & 0 & 0& 0 &0 & 0 & $x_5$ & 0 & 0\\ \hline
$\beta_*^{r_2}$ & 0 & $a_2$ &0 & $-$&0&$-$ &0&$b_{10}$ \\ \hline
$\beta_*^{r_2+1}$ &$-$&$-$&0 &$-$  & $a_4$ &$-$ &$-$& $-$\\ \hline
\end{tabular}
\end{center}
where $r_2=\nu_2(q^2-1)$. 
\end{itemize}
\end{Th}

Next we have the analogous calculation for $B\Sol(q)$. 

\begin{Th}\label{Sol(q)}
Fix an odd prime power q. Then 
\[H_*(\Omega B\Sol(q)) \cong P[a_6]/(a_6^2) \otimes P[b_{10}, c_{12}, e_{26}]\otimes E[y_7, y_{11}, y_{13}]\otimes P[y_{14}]/(y_{14}^2),\]
as a module over $H_*(\Omega DI(4))\cong P[a_6]/(a_6^2) \otimes P[b_{10}, c_{12}, e_{26}]$. 
Furthermore:
\begin{itemize}
\item The relations  which determine the algebra extension are given by $y_7^2 = y_{11}^2 = y_{13}^2 = y_{14}^2 = 0$, $[a_6, y_{14}] = b_{10}^2$, $[b_{10},y_{14}] = c_{12}^2$, and $[c_{12}, y_{14}] = e_{26} + a_6b_{10}^2$, and $[e_{26}, y_{14}] = b_{10}^4$. All other commutators of generators are trivial. 
\item The reduced coproduct is given by $\widebar{\Delta}(c_{12}) = a_6\otimes a_6$, and  $\widebar{\Delta}(y_{14}) = y_7\otimes y_7$. All other generators are primitive.
\item The action of the dual Steenrod algebra is determined by 
\begin{center}
 \begin{tabular}{||l|c|c|c|c|c|l|c|l||} \hline
 & $a_6$ & $y_7$ & $b_{10}$ & $y_{11}$ &   $c_{12}$ &   $y_{13}$ & $y_{14}$ & $e_{26}$   \\ \hline
$Sq^1_*$ & $0$ & $0$ & $0$ & $0$ &   $0$ & $0$  & $y_{13}$ & $0$   \\ \hline
$Sq^2_*$ & $0$ & $0$ & $0$ & $0$  & $b_{10}$ &  $y_{11}$ & $0$ & $c_{12}^2$  \\ \hline
$Sq^4_*$ & $0$ & $0$& $a_6$ & $y_7$  &  $0$ &   $0$ & $0$ & $0$  \\ \hline
\end{tabular}
\end{center}
and the Steenrod axioms.
\item The homology Bockstein spectral sequence is determined by the table, 
\begin{center}
 \begin{tabular}{||l|c|c|c|c|c|c|c|c|l||} \hline
& $a_6$ & $y_7$ & $b_{10}$ & $y_{11}$  & $c_{12}$ &$y_{13}$& $y_{14}$ & $a_6y_7$ & $y_{13}y_{14}$ \\\hline
$Sq_*^1$ &$0$&$0$&$0$&$0$&$0$&$0$&$y_{13}$&$0$&$0$ \\ \hline
$\beta_*^{r_4-1}$ &$0$&$0$&$0$&$b_{10}$&$0$&$-$&$-$&$0$&$e_{26}$\\ \hline
$\beta_*^{r_4}$ &$0$&$a_6$&$-$&$-$&$0$&$-$&$-$&$0$&$-$\\ \hline
$\beta_*^{r_4+1}$ &$-$&$-$&$-$&$-$&$0$&$-$&$-$&$c_{12}$&$-$  \\ \hline
\end{tabular}
\end{center}
where $r_4 = \nu_2(q^4-1)$.
\end{itemize}
\end{Th}
  
The paper is organized as follows. In Section 1 we record some basic facts which are the basis for our calculation. The loop space homology of $BG_2(q)\doscom$ and $B\Sol(q)$ are calculated in Sections 2 and 3 respectively.

Some of the calculations presented here can be carried out more easily using the general methods developed recently by Daisuke Kishimoto and Akira Kono \cite{KK}. The authors are very grateful to Kono for the interest he showed in our results and for pointing out an error in the calculation of the algebra structures in an earlier version of this paper.


\section{Preliminaries}
Recall the Quillen-Friedlander fibre square \cite{Fr} for groups of Lie type. If $G$ is a complex reductive Lie group, and $G(q)$ is the corresponding algebraic group over the field of $q$ elements, then after completion at a prime $p$ not dividing $q$, there is a homotopy fibre square
\begin{equation}\label{FFS}
\begin{diagram}
BG(\F_q)\pcom && \rTo && BG\pcom\\
\dTo &&&& \dTo_{\Delta}\\
BG\pcom &&\rTo{1\top\psi^q}&& BG\pcom\times BG\pcom
\end{diagram}
\end{equation}
where $\psi^q$ is the $q$-th unstable Adams operation, and $\Delta$ is the diagonal map.  In particular, since for any self map $f\colon X\to X$, $\mathrm{hofib}(X\rTo{1\top f} X\times X) \simeq \Omega X$, one has a fibration sequence of loop spaces and loop maps:
\begin{equation}\label{selfmap-fib}
\Omega BG(q)\pcom\rTo G\pcom\rTo{f^q} G\pcom.
\end{equation}
All $p$-compact groups, in particular $DI(4)$, admit unstable Adams operations of degree $q$, where $q$ is a $p$-adic unit. The corresponding fibre square for $DI(4)$ was used by Benson to define $B\Sol(q)$ \cite{Benson}.

We next record three well known cohomology algebras, which will be used in our calculation. As a convention we will use Roman alphabet to denote classes in mod 2 homology and cohomology, and Greek letters to denote classes in integral homology and cohomology. A good reference for the cohomology of Lie groups is \cite{Mimura-Toda}. 
\bigskip

\noindent{\bf The Spaces $\mathbf{BSU(3)}$, $\mathbf{SU(3)}$ and $\mathbf{\Omega SU(3)}$.}
\begin{equation}\label{BSU(3)}
H^*(BSU(3))\cong P[u_4, u_6]\quad\mathrm{and}\quad H^*(BSU(3),\Z) \cong P[\gamma_4,\gamma_6],\end{equation}
both as algebras, with $Sq^2(u_4) = u_6$. Recall also that 
\begin{equation}\label{SU3}
H_*(SU(3),\Z)\cong E[\chi_3,\chi_5]\end{equation}
 as a Hopf algebra. An elementary calculation, using the EMSS, yields
\begin{equation}\label{loopsSU3}
H_*(\Omega SU(3),\Z)\cong P[\alpha_2,\alpha_4]\end{equation}
 as an algebra, with the Hopf algebra structure determined by 
$\widebar{\Delta}(\alpha_4)=\alpha_2\otimes \alpha_2$, where $\widebar{\Delta}$ denotes the reduced diagonal. Since these
algebras are torsion free, 
\[H_*(SU(3))\cong H_*(SU(3),\Z)\otimes
\F_2\cong E[x_3,x_5],\] and \[H_*(\Omega SU(3)) \cong H_*(\Omega SU(3),\Z)\otimes \F_2\cong P[a_2,a_4]\]
as Hopf algebras, with $Sq_*^2(x_5)=x_3$ and $Sq_*^2(a_4)=a_2$. 

\bigskip

\noindent{\bf The Spaces $\mathbf{BG_2}$ and $\mathbf{G_2}$.}
\begin{equation}\label{BG2}
H^*(BG_2) \cong P[u_4,u_6, t_7],\end{equation}
with $Sq^2(u_4) = u_6$ and $Sq^1(u_6) = t_7$. These are the rank 3 mod 2 Dickson invariants. The group $SU(3)$ is a subgroup of $G_2$ and the inclusion induces the obvious projection on mod 2 cohomology. Recall also
that 
\begin{equation}\label{G2-hlgy}
H_*(G_2)\cong E[x_3,x_5,x_6]\end{equation}
 with $Sq_*^1(x_6)=x_5$,
$Sq_*^2(x_5)=x_3$ and $\widebar{\Delta}(x_6) = x_3\otimes x_3$.

\bigskip

\noindent{\bf The Spaces $\mathbf{BDI(4)}$ and $\mathbf{DI(4)}$.}
Let $BDI(4)$ denote the classifying space of the 2-compact group $DI(4)$ \cite{DW}. Thus 
\begin{equation}\label{BDI4}
H^*(BDI(4), \F_2) \cong P[v_8, v_{12}, v_{14}, s_{15}],\end{equation}
with $Sq^4(v_8)=v_{12}$, $Sq^2(v_{12}) = v_{14}$, and $Sq^1(v_{14}) = s_{15}$. 
One also has
\begin{equation}\label{DI4}
H_*(DI(4),\F_2) = E[y_7,y_{11},y_{13},y_{14}]\end{equation}
with $Sq^4(y_7)=y_{11}$, $Sq^2(y_{11}) = y_{13}$,  $Sq^1(y_{13}) = y_{14}$, and $\widebar{\Delta}(y_{14}) = y_7\otimes y_7$.


\section{Loop space homology of $BG_2(q)\doscom$.}
Next we calculate the loop space homology of $BG_2(q)\doscom$. To
avoid an awkward notation, we use $\widehat{(-)}$ to denote $(-)\doscom$ where it makes sense to do so.  

The mod-2 loop space homology of  $G_2$, which is necessary for the calculation in hand, is well known (see for instance \cite{Bott, Kono-Kozima}), but we include a brief calculation here for the convenience of the reader. 

There is a
fibration 
\begin{equation}\label{e2}
SU(3)\rTo G_2\rTo S^6.
\end{equation}
To calculate the loop space homology of $G_2$, we use the mod-2 and
integral
homology Serre spectral sequences for the
fibration obtained from (\ref{e2}) by looping. Consider first the integral spectral sequence
\[E^2 = H_*(\Omega SU(3),\Z)\otimes H_*(\Omega S^6,\Z) \cong
P[\alpha_2,\alpha_4]\otimes T[\beta_5].\]
Differentials in this spectral sequence respect the coproduct structure, and hence $d_5(\beta_5)$ must be
a primitive element in $H_4(\Omega SU(3),\Z)$.  One has
$\bar{\Delta}(\alpha_4)=\alpha_2\otimes \alpha_2$, and $\bar\Delta(\alpha_2^2) =
2\alpha_2\otimes \alpha_2$, hence any primitive in $H_4(\Omega SU(3), \Z)$ is a multiple of $2\alpha_4-\alpha_2^2$, and so $d_5(\beta_5) = A(2\alpha_4-\alpha_2^2)$, for some $A\in\Z$.

To determine the value of $A$, consider 
the fibration
\[\Omega S^6\rTo{\delta} SU(3) \rTo G_2.\]
An easy calculation with the homology SSS of this fibration, using the fact that $H_5(G_2,\Z)=\Z/2$,  shows
that $\delta_*(\beta_5) = 2x_5$. Hence using the commutative diagram
\[\begin{diagram}\Omega SU(3) & \rTo & \Omega G_2 & \rTo &  \Omega S^6\\
\dTo{=} && \dTo && \dTo{\delta}\\
\Omega SU(3) & \rTo & * & \rTo &  SU(3)
\end{diagram}\]
and naturality of the SSS,  one has in the top fibration
$d_5(\beta_5) = 2\alpha_4$ modulo the ideal generated by $\alpha_2$. This shows that
$A=1$.

Reducing this calculation mod 2, one has 
\[E^2 = H_*(\Omega SU(3))\otimes H_*(\Omega S^6) \cong
P[a_2,a_4]\otimes T[b_5],\]
and it follows easily that $d_5(b_5) = a_2^2$.  
 Thus one has
\begin{equation}\label{loopsG2}
H_*(\Omega G_2) \cong P[a_2, a_4]/(a_2^2)\otimes P[b_{10}]\end{equation}
as a module over $H_*(\Omega SU(3))$. Since the element $b_{10}$ has infinite height already in the $E^\infty$ page of the spectral sequence, it is also an element of infinite height in $H^*(\Omega G_2)$, and the structure given in (\ref{loopsG2}) is the algebra structure. Notice that the element $b_{10}$ is determined only up to an additive summand of an element in the image  of restriction from $H_*(\Omega SU(3))$. (Notice also that in integral homology one has the relation $a_2^2 = 2a_4$.)

The Hopf algebra structure of $H_*(\Omega G_2)$ is determined by $\bar\Delta(a_4) = a_2\otimes a_2$, while $b_{10}$ can be chosen to be primitive (if some choice of $b_{10}$ is not primitive, then $b_{10}' = b_{10}+a_4^2a_2$ is).  It follows that in cohomology the dual of $a_4$ is the square of the dual of $a_2$, and so one has $Sq^2_*(a_4) = a_2$.

Dually, one has 
\[H^*(\Omega G_2) \cong P[\widebar{a}_2]/(\widebar{a}_2^4)\otimes \Gamma[a_8, \widebar{b}_{10}],\]
where $\widebar{a}_2$ and $\widebar{b}_{10}$ are the duals of $a_2$ and $b_{10}$ respectively, and $a_8$ is dual to $a_4^2$.

Deciding the action of the homology Steenrod squares on $b_{10}$ requires more calculation. The authors are grateful to Akira Kono for sketching for them the argument that follows. Let $\widetilde{G_2}$ denote the 3-connected cover of $G_2$. Thus there is a principal fibration
\[K(\Z,2)\rTo \widetilde{G_2} \rTo G_2.\]
Using the mod-2 cohomology SSS  for this fibration,  an elementary computation shows that 
\[H^*(\widetilde{G_2}) \cong P[u_8] \otimes E[y_9, z_{11}],\]
where $u_8$ restricts to $\iota_2^4\in H^8(K(\Z,2))$. The classes $y_9$ and $z_{11}$ correspond to the infinite cycles in the spectral sequence given by $\iota_2^2b_5$ and $\iota_2 a_3^3$, where $a_3$ and $b_5$ denote the generators of $H^*(G_2)$. By analysing the SSS for the fibration
\[\widetilde{G_2} \rTo G_2 \rTo K(\Z, 3),\] it is easy to see that $Sq^1(u_8) = y_9$ and $Sq^2(y_9) = z_{11}$. Finally, using the spectral sequence for 
\[\Omega G_2\rTo K(\Z,2) \rTo \widetilde{G_2}\]
one observes that $\iota_2$ restricts to $\widebar{a}_2$, and so $\iota_2^4$ restricts trivially, and is therefore the image of $u_8$ under the inflation map. The rest of the spectral sequence is determined by letting $y_9$ and $z_{11}$ be the image of $\widebar{a}_8$ and $\widebar{b}_{10}$ under the  transgression. In particular  it follows that $Sq^2(\widebar{a}_8) = \widebar{b}_{10}$. Dually in homology we have $Sq_*^2(b_{10}) = a_4^2$. 

To calculate the loop space homology of $BG_2(q)\doscom$, consider first the Friedlander fibre square (\ref{FFS}) for $G_2$. Taking iterated fibres on the left column of the square we get a sequence of fibrations
\[\Omega\widehat{G_2} \rTo \Omega BG_2(q)\doscom \rTo \widehat{G_2} \rTo{f^q} \widehat{G_2}\rTo BG_2(q)\doscom\rTo B\widehat{G_2},\]
where $\widehat{G_2}$ denotes $(G_2)\doscom$.
For the actual calculation, we use the SSS for the fibration
\begin{equation}\label{loops-bg2q-fib}\Omega \widehat{G_2} \rTo \Omega BG_2(q)\doscom \rTo \widehat{G_2}.\end{equation}
Thus, we start by calculating the map induced by $f^q$ on homology.

Expanding the Friendlander fibre square, one sees that $f^q$ is the composite
\[\widehat{G_2} \rTo{1\top\Omega\psi^q} \widehat{G_2}\times \widehat{G_2}\rTo{-1\top 1} \widehat{G_2}\times \widehat{G_2}\rTo{\mu} \widehat{G_2}.\]
There is an isomorphism of modules over $P[u_4, u_7, u_6^2]$,
\[H^*(BG_2)\cong P[u_4, u_7, u_6^2]\otimes E[u_6],\]
with $Sq^1(u_6) = u_7$. Considering this as a differential graded algebra with the differential given by $Sq^1$, and taking cohomology, we get the $E_2$ term of the BSS for $H^*(BG_2)$,
\[E_2  \cong P[u_4, u_6^2],\]
which is concentrated in even degrees, and so  $E_2=E_\infty$. Hence the integral cohomology of $BG_2$ is given by 
\[H^*(BG_2, \Z)\cong P[u_4, u_7, v_{12}]/(2u_7).\]
Notice that $u_4$ and $v_{12}$ are torsion free classes and $(\psi^q)^*(u_4) = q^2u_4$, while $(\psi^q)^*(v_{12}) = q^6v_{12}$. On the other hand, since the class $u_7$ is of order 2, every element in the ideal it generates is of order 2, and since $\psi^q$ is a mod-2 equivalence for $q$ odd, the ideal generated by $u_7$ is fixed under $(\psi^q)^*$.

Using the BSS for $H_*(G_2)$ we see that 
\[H_i(G_2, \Z) \cong \begin{cases}
								\Z & i=0, 3, 11, 14\\
								\Z/2 & i = 5, 8
					\end{cases},\]
					while all other homology groups vanish.
Let $\chi_i$ denote a generator for $H_i(G_2,\Z)$ for those values of $i$, where the respective homology group is nontrivial.  Using the known algebra structure of $H_*(G_2)$, it is easy to conclude that $\chi_3$, $\chi_5$ and $\chi_{11}$ are indecomposable, and that $\chi_3\chi_5 = \chi_8$ and $\chi_3\chi_{11} = \chi_{14}$. 
Using the SSS for the path loop fibration over $BG_2$ and naturality, we conclude that $\Omega\psi^q_*(\chi_3) = q^2\chi_3$ and $\Omega\psi^q_*(\chi_{11}) = q^6\chi_{11}$, while $\Omega\psi^q_*(\chi_5) = \chi_5$ and $\Omega\psi^q_*(\chi_8) = \chi_8$. In mod-2 homology both $(\psi^q)_*$ and $(\Omega\psi^q)_*$ behave like the identity.

Using the K\"unneth formula, we see that for $n\le 11$
\[H_n(G_2\times G_2,\Z) \cong \bigoplus_{i+j=n}H_i(G_2,\Z)\otimes H_j(G_2,\Z).\]

This, and the information about $\Omega\psi_*^q$ allows us to easily calculate $f^q_*$ on $H_*(G_2, Z)$. One has 
\[f^q_*(\chi_3) = (q^2-1)\chi_3,\quad f^q_*(\chi_5) = 0,\quad \mathrm{and}\quad f^q_*(\chi_{11}) = (q^6-1)\chi_{11}.\]
On mod-2 homology $f^q_*$ is trivial.

Now consider the fibration (\ref{loops-bg2q-fib}), which is induced from the path-loop fibration over $G_2$ via the map $f_q$. 
Since $f^q_*$ is trivial on mod-2 homology, the SSS for (\ref{loops-bg2q-fib}) collapses at $E^2$, and  it follows that  for all odd $q$ there is an isomorphism of modules over $H_*(\Omega G_2)$:
\begin{equation}\label{g2qhlgy-2}H_*(\Omega BG_2(q)\doscom) \cong \left\{P[a_2]/(a_2^2)\otimes
P[a_4,b_{10}]\right\}\otimes E[x_3,x_5]\otimes P[z_6]/(z_6^2).\end{equation}

The structure of $H_*(\Omega BG_2(q)\doscom)$ as a module over the dual Steenrod algebra follows from the information we have about the two factors. Namely, $Sq_*^2(a_4) = a_2$, $Sq^2_*(x_5) = x_3$, $Sq^2_*(b_{10}) = a_4^2$, and $Sq^1_*(z_6) = x_5$. This is summarised in the following table.
\begin{center}
 \begin{tabular}{||l|c|c|c|l||} \hline
$q\equiv 1(4)$ & $a_4$ & $x_5$  & $z_6$ & $b_{10}$ \\ \hline
$Sq^1_*$  & 0 &0 &  $x_5$ & 0 \\ \hline
$Sq^2_*$  & $a_2$ & $x_3$ & 0 & $a_4^2$ \\ \hline
\end{tabular}\end{center}

The classes $a_2$ and $x_3$ are primitive with respect to the diagonal in $H_*(\Omega BG_2(q)\doscom)$ for dimensional reasons. The class $x_5$ can be chosen to be primitive, since for any choice $\widebar{\Delta}(x_5) = A(a_2\otimes x_3 + x_3\otimes a_2)$, for some $A\in \F_2$, then $x_5 + Aa_2x_3$ is primitive, has the same action of $Sq_*^2$ as $x_5$, and represents the same class modulo $H_*(\Omega G_2)$. One also has $\widebar{\Delta}(a_4) = a_2\otimes a_2$, since $H_*(\Omega G_2)$  is a Hopf subalgebra. Finally, the class $z_6$ can be chosen to have reduced diagonal $x_3\otimes x_3$, since for any choice of representative, one has $\widebar{\Delta}(z_6) = x_3\otimes x_3 + B(a_2\otimes a_4 + a_4\otimes a_2)$, and so $z'_6 = z_6 + Ba_2a_4$ is congruent to $z_6$ modulo $H_*(\Omega G_2)$, has the same action of $Sq_*^1$ as $z_6$, and has the required diagonal.

Next, we compute the algebra extension. To do that, fix the representatives for $x_3, x_5$ and $z_6$ as above. Notice first that $x_3^2=0$ and $z_6^2=0$, since there are no primitives in the respective dimensions. Similarly, $x_5$ is primitive, and so $x_5^2 = Ab_{10}$ for some $A\in \F_2$. Applying $Sq_*^2$ to both sides we see that $A=0$, so $x_5^2=0$.

Next we  systematically examine all the commutators involving $x_3, x_5$ and $z_6$, as listed in the following table.
\begin{center}
 \begin{tabular}{||l|c|c|c|c|c|c|c|l||} \hline
  5 & 7 & 8 & 9 & 10 & 11 & 13 & 15 & 16\\ \hline 
  $[a_2,x_3]$ & $[a_2, x_5]$ & $[a_2, z_6]$&$[a_4, x_5]$& $[a_4, z_6]$ &$[x_5, z_6]$ & $[x_3, b_{10}]$& $[x_5,b_{10}]$& $[z_6, b_{10}]$ \\
                      & $[x_3, a_4]$ & $[x_3, x_5]$ & $[x_3,z_6]$&&&&& \\\hline                 
 \end{tabular}
 \end{center}
 We will show that 
 \[[a_2,z_6] = a_4^2,\quad [a_4,z_6]=b_{10} + a_4^2a_2,\quad\text{and} \quad [b_{10},z_6] = a_4^4,\] while all the other commutators in the table vanish. 
 
Observe first that every non-primitive class among the generators of $H_*(\Omega BG_2(q)\doscom)$ has a reduced diagonal consisting of a single element. The commutator of two primitives is always primitive, and if $a$ is a primitive and $\widebar{\Delta}(b) = c\otimes c$, then 
\begin{equation}\label{com-for}\widebar{\Delta}([a,b]) = c\otimes [a,c]  + [a,c]\otimes c.\end{equation}
Thus the commutators $[a_2,x_3]$, $[a_2,x_5]$, $[x_3,x_5]$, $[x_3, b_{10}]$ and $[x_5, b_{10}]$ are automatically primitive. 

Since $[a_2, x_3] = Ax_5$. Applying $Sq_*^2$ to both sides, it follows that $A=0$, so $a_2$ and $x_3$ commute. Since $\widebar{\Delta}(a_4)  = a_2\otimes a_2$ and $\widebar{\Delta}(z_6) = x_3\otimes x_3$,  (\ref{com-for}) applies, and since $[a_2,x_3]=0$, $[x_3,a_4]$ and $[a_2,z_6]$ are also primitive. The only other primitive in dimension 7 is $[a_2, x_5]$, and so $[x_3, a_4] = A[a_2, x_5]$, for some $A\in\F_2$. Similarly 
\[[a_2,z_6] = B[x_3,x_5] + Ca_4^2,\]
as $[x_3,x_5]$ and $a_4^2$ are the only other primitives in dimension 8. But $Sq^1_*$ applied to both sides yields 0 on the right hand side, and $[a_2, x_5]$ on the left hand side. Hence $a_2$ commutes with $x_5$, and consequently $a_4$ commutes with $x_3$. 
 
By a similar method we analyse $[a_4,z_6]$, $[a_4, x_5]$, $[x_3, z_6]$, $[x_3, x_5]$ and $[x_5,z_6]$. First, by direct calculation, and the results already listed above, 
\[\widebar{\Delta}([a_4, z_6]) = a_2\otimes[a_2,z_6] + [a_2,z_6]\otimes a_2.\]
Thus before we have decided whether $a_2$ commutes with $z_6$, we must assume $[a_4,z_6] = Ab_{10} + Ba_4^2a_2$, for some $A, B\in \F_2$. Applying $Sq^1_*$ to both sides we have $[a_4, x_5]=0$. 
Now, since $x_3$ and $x_5$ commute, $[x_5, z_6]$ is primitive. Since $x_5$ commutes with both $a_2$ and $a_4$, there are no other nonzero primitives in dimension 11, and so $[x_5,z_6]=0$. Applying $Sq_*^2$ and then $Sq_*^1$, we conclude that $[x_3,z_6]$ and $[x_3,x_5]$ both vanish.

Next, notice that $[z_6, b_{10}]$ is a primitive class. Hence $[z_6,b_{10}] = Aa_4^4$ for some $A\in\F_2$. Applying $Sq_*^1$ and then $Sq_*^2$ to both sides, and using the fact that $x_5$ and $a_4$ commute, we conclude that $b_{10}$ commutes with $x_5$ and $x_3$. 

It remains to analyse the commutators $[a_2,z_6]$, $[a_4, z_6]$ and $[b_{10},z_6]$. To do that, recall from \cite{Grbic},
\[H^*(BG_2(q)) = P[u_4, u_6, t_7, y_3, y_5]/(y_5^2 + y_3t_7 + y_3^2u_4,\, y_3^4 +y_5t_7 + y_3^2 u_6).\]
Denote classes in $H_*(BG_2(q))$ by adding a bar to the corresponding cohomology class, and consider the cobar spectral sequence for $H_*(\Omega BG_2(q)\doscom)$. Thus 
\[E^2 \cong \Cotor^{H_*(\Omega BG_2(q)\doscom)}(\F_2,\F_2) \cong H_*(T(\Sigma^{-1}(\widebar{H}_*(BG_2(q)\doscom)), d_E),\]
where $d_E$ is the differential on $T(\Sigma^{-1}(\widebar{H}_*(BG_2(q)\doscom))$ induced by the reduced diagonal. If $x, y\in \widebar{H}_*(BG_2(q)\doscom)$ are any classes, we denote the corresponding elements of the tensor algebra by $[x]$ $[y]$ etc., and their product in the tensor algebra structure by standard bar notation $[x|y]$. 

Consider the homology classes $\widebar{y_5^2}$, $\widebar{y_3t_7}$ and  $\widebar{y_3^2u_4}$. The corresponding reduced diagonals are $\widebar{y_5}\otimes\widebar{y_5}$, $\widebar{y_3}\otimes\widebar{t_7} + \widebar{t_7}\otimes\widebar{y_3}$ and $\widebar{y_3}\otimes \widebar{y_3u_4} + \widebar{y_3u_4}\otimes \widebar{y_3} + \widebar{y_3^2}\otimes \widebar{u_4} + \widebar{u_4}\otimes \widebar{y_3^2}$ respectively. Hence
\[d_E([\widebar{y_5^2}]) = [\widebar{y_5}]^2, \quad d_E([\widebar{y_3t_7}]) = [[\widebar{y_3}],[\widebar{t_7}]], \quad\mathrm{and}\quad d_E([\widebar{y_3^2u_4}]) = [[\widebar{y_3}],[\widebar{y_3u_4}]] + [[\widebar{y_3^2}],[\widebar{u_4}]].\]
On the other hand, since $y_5^2 + y_3t_7 + y_3^2u_4=0$ in $H^*(BG_2(q))$, we have 
\[[\widebar{y_5}]^2 + [[\widebar{y_3}],[\widebar{t_7}]] =  [[\widebar{y_3}],[\widebar{y_3u_4}]] + [[\widebar{y_3^2}],[\widebar{u_4}]].\] Furthermore, $[\widebar{y_5}] $ and $[\widebar{y_3}]$ and $[\widebar{t_7}]$ are all cycles, which are permanent for dimensional reasons and hence represent  $a_4$, $a_2$ and $z_6$ respectively in loop space homology. The equations above show that the expression $[\widebar{y_5}]^2 + [[\widebar{y_3}],[\widebar{t_7}]]$ is a boundary, and so we obtain the relation $[a_2,z_6] = a^2_4$. Next, notice that $Sq_*^2([a_4,z_6]) = [a_2,z_6] = a_4^2$, while $\widebar{\Delta}([a_4, z_6]) = a_2\otimes a_4^2 + a_4^2\otimes a_2$. Hence we conclude that 
\[[a_4,z_6] = b_{10} + a_4^2a_2 = b_{10} + [a_2,z_6]a_2.\]
Finally,  since $b_{10} =[a_4,z_6]+a_4^2a_2$, we calculate directly, 
\begin{multline*}[b_{10}, z_6] = [[a_4,z_6]+a_4^2a_2, z_6]=[[a_4,z_6], z_6] + [a_4^2a_2, z_6] =\\ [a_4,z_6^2] + [[a_2,z_6]a_2,z_6] = 0 + [a_2,z_6]^2 = a_4^4.
\end{multline*}
This completes the computation of the Hopf algebra structure. To summarise, we have shown that 
\[H_*(\Omega BG_2(q)\doscom)\cong P[a_2]/(a_2^2)\otimes P[a_4, b_{10}]\otimes E[x_3,x_5]\otimes P[z_6]/(z_6^2),\]
as modules over $H_*(\Omega BG_2)$. The relations which determine the algebra extension are given by $x_3^2 = x_5^2 = z_6^2 =0$, $[a_2,z_6] = a_4^2$, $[a_4,z_6] = b_{10} + a_2a_4^2$, and $[b_{10}, z_6] = a^4_4$. All other commutators of generators are trivial. The coproduct is given by $\widebar{\Delta}(a_4) = a_2\otimes a_2$, $\widebar{\Delta}(z_6) = x_3\otimes x_3$, and all other generators are primitive.

It remains to compute the Bockstein spectral sequence for $H_*(\Omega BG_2(q)\doscom)$. This is done by calculating the integral SSS for the fibration in the top row of the diagarm

\[\begin{diagram}
\Omega \widehat{G}_2 & \rTo & \Omega BG_2(q)\doscom & \rTo &\widehat{G}_2\\
\dIgual && \dTo && \dTo{f^q}\\
\Omega \widehat{G}_2 & \rTo & * & \rTo &\widehat{G}_2
\end{diagram}\]
We use naturality and the action of $f^q_*$ computed above. First, analyse the SSS for the bottom row, using the same notation we have been using before. One has $d_3(\chi_3) = a_2$, and since in integral homology $a_2^2 = 2a_4$, $d_3(a_2\chi_3) = 2a_4$. Hence $d_3(a_4^k\chi_3) = a_4^ka_2$, and $d_3(a_4^ka_2\chi_3) = 2a_4^{k+1}$. Since $\chi_8 = \chi_3\chi_5$, it follows that $d_3(\chi_8) = a_2\chi_5$, but $d_3(a_2\chi_8) = 0$. Similarly, $d_3(\chi_{14}) = a_2\chi_{11}$. In addition one must have $d_3(\chi_{11}) = a_2\chi_8$, since otherwise $a_2\chi_8$ will be an infinite cycle. This determines $d_3$. The next nontrivial differential is $d_5$, which takes $\chi_5$ isomorphically to $a_4$, which in $E^5$ is a class of order 2. Finally $d_{11}$ takes $\chi_{11}$ to $b_{10}$, and $E^{12} = E^\infty$. Now, using naturality of the spectral sequence and our knowledge of $f^q_*$, it follows that in the SSS for the top row in the diagram, $d_3(\chi_3) = (q^2-1)a_2$, $d_5(\chi_5) = 0$, and $d_{11}(\chi_{11}) = (q^6-1)b_{10}$. This information suffices for the computation of the BSS. The integral calculation yields in particular the observation that $a_2$ is a class of order $(q^2-1)$, while $b_{10}$ has order $(q^6-1)$.

Now, consider $H_*(\Omega BG_2(q)\doscom)$ as a module over $H_*(\Omega G_2)$ as in (\ref{g2qhlgy-2}):
\[H_*(\Omega BG_2(q)\doscom) \cong \left\{P[a_2]/(a_2^2)\otimes
P[a_4,b_{10}]\right\}\otimes E[x_3,x_5]\otimes P[z_6]/(z_6^2).\]
Taking the $Sq_*^1$ homology, one has 
\[E^2 \cong 
P[a_2]/(a_2^2)\otimes P[a_4, b_{10}]\otimes E[x_3, h_{11}].\]
where the class $h_{11}$ is represented by the $Sq_*^1$ cycle $x_5z_6$. Notice that for any odd $q$, $r_6=\nu_2(q^6-1) = \nu_2(q^2-1) = r_2$. Hence the next nontrivial Bockstein operator is $\beta_*^{r_2}(x_3) = a_2$, and $\beta^{r_2}_*(h_{11}) = b_{10}$. It now follows that $\beta^{r_2+1}_*(a_2x_3) = a_4$, and $E^{r_2+2} = E^\infty$.

The results are summarised in the following table.
\begin{center}
 \begin{tabular}{||l|c|c|c|c|c|c|c|l||} \hline
$q\equiv 1(4)$& $a_2$ & $x_3$ & $a_4$ & $x_5$  & $a_2x_3$ & $z_6$ & $b_{10}$ & $x_5z_6$\\ \hline
$Sq^1_*$  & 0 & 0& 0 &0 & 0 & $x_5$ & 0 & 0\\ \hline
$\beta_*^{r_2}$ & 0 & $a_2$ &0 & $-$&0&$-$ &0&$b_{10}$ \\ \hline
$\beta_*^{r_2+1}$ &$-$&$-$&0 &$-$  & $a_4$ &$-$ &$-$& $-$\\ \hline
\end{tabular}
\end{center}

This completes the proof of Theorem \ref{G2}.


\section{Loop space homology of $B\Sol(q)$}

For any odd prime power $q$, the 2-local finite group $\Sol(q)$ is defined in \cite{LO}. The starting point is a family of saturated fusion systems $\calf_{\Sol(q)}$ over the Sylow 2-subgroup of $\mathrm{Spin}_7(q)$. These fusion systems were originally defined by Ron Solomon \cite{Sol} as a part of his contribution to the classification of finite simple groups. He did not use the language of fusion systems, but essentially presented the entire family and studied its general behaviour. In \cite{Benson}, Benson gave a construction of a family of spaces $B\Sol(q)$, which he claimed realise the fusion patterns defined by Solomon. These spaces are given as the pullback spaces in the diagram
\begin{equation}\label{fib-square}
\begin{diagram}B\Sol(q) && \rTo & &BDI(4) \\
\dTo &&&& \dTo_{\Delta}\\
BDI(4)&& \rTo{\psi^q\top 1} && BDI(4)^{\times 2}
\end{diagram}\end{equation}
 where $\Delta$ is the diagonal and $\psi^q$ is the degree $q$ unstable Adams operation on $BDI(4)$ constructed by Notbohm \cite{Notbohm}. 
The paper \cite{LO} unifies the two constructions. On one hand it is shown that the fusion patterns defined by Solomon are indeed saturated fusion systems, each of which admits an associated centric linking system $\call_{Sol(q)}$, and on the other hand that the classifying spaces of the corresponding $2$-local finite groups $B\Sol(q)\defeq |\call_{\Sol(q)}|\doscom$ coincide with the spaces  constructed by Benson, whose approach allows a calculation of the mod 2 cohomology of $B\Sol(q)$, as demonstrated in \cite{Grbic}. We shall also utilize Benson's pullback diagram in the current work. 

In what follows we will denote $\Omega BDI(4)$ by $DI(4)$ (not to be confused with the notation $DI(n)$ which is sometimes used to denote the rank $n$ algebra of Dickson invariants).

As before, one has a fibration of loop spaces and loop maps
\[\Omega DI(4) \rTo \Omega B\Sol(q) \rTo DI(4),\]
resulting from looping the left the left hand side column in Benson's pull back diagram (\ref{fib-square}).
Thus our first task is to compute the loop space homology of $DI(4)$. 

\begin{Prop}\label{Omega-DI(4)-homology}
There is an isomorphism of Hopf algebras
\[H_*(\Omega DI(4)) \cong P[a_6]/(a_6^2) \otimes P[b_{10}, c_{12}, e_{26}],\]
where $a_6$, $b_{10}$ and $e_{26}$ are primitive, and $\widebar{\Delta}(c_{12}) = a_6\otimes a_6$. The action of the dual Steenrod algebra is determined by $Sq_*^4(b_{10}) = a_6$, $Sq_*^2(c_{12}) = b_{10}$, and $Sq^2_*(e_{26}) = c_{12}^2$.
\end{Prop}
\begin{proof}
Consider the homology EMSS for the path-loop fibration over $DI(4)$. The $E^2$ term is given by 
\[\Cotor^{H_*(DI(4))}(\F_2,\F_2) \cong \Ext_{H^*(DI(4))}(\F_2,\F_2).\]
The isomorphism holds since $H_*(DI(4))$ is of finite type. To calculate the right hand side, consider the differential graded Hopf algebra
\[ P_* \defeq (P[x_7]/(x_7^4)\otimes E[y_{11}, z_{13}])\otimes (P[\hat{a}_6]/(\hat{a}^4)\otimes \Gamma[\hat{b}_{10}, \hat{t}_{24}, \hat{e}_{26}]),\]
where the left factor is primitively generated, and the right factor is the dual of the Hopf algebra $P[a_6]/(a^2)\otimes P[b_{10}, c_{12}, e_{26}]$, where all generators but $c_{12}$ are primitive, and $\widebar{\Delta}(c_{12}) = a_6\otimes a_6$.
We denote by $\gamma_k(\hat{b})$ the generator of $\Gamma[\hat{b}_{10}]$ in dimension $10k$, where $\gamma_1(\hat{b}) = b_{10}$, and $\gamma_0(\hat{b})=1$. We use similar notation for the generators corresponding to $\hat{c}_{12}$ and $\hat{e}_{26}$. Thus as an $\F_2$-algebra $\Gamma[\hat{b}_{10}, \hat{t}_{24}, \hat{e}_{26}]$ can be written as 
\[\bigotimes_{n\geq 0}E[\gamma_{2^n}(\hat{b}), \gamma_{2^n}(\hat{t}), \gamma_{2^n}(\hat{e})],\]
where we omit subscripts for short.
The differential on $P_*$ is given on generators by 
\begin{itemize}
\item $d(x) = d(y) = d(z) = 0$, 
\item $d(\hat{a}) = x$,
\item $d(\gamma_{2^n}(\hat{b})) = y\gamma_{2^n-1}(\hat{b})\;$, 
\item $d(\hat{a}^2) = z$,
\item $d(\gamma_{2^n}(\hat{t})) = z\hat{a}^2 \gamma_{2^n-1}(\hat{c})\;$, and 
\item $d(\gamma_{2^n}(\hat{e})) =  x^3  \hat{a}\gamma_{2^n-1}(\hat{e})$,
\end{itemize}
and is required to satisfy the Leibniz rule on products. Notice that the differential is, in particular, a map of graded algebras over  $H^*(DI(4)) = P[x_7]/(x_7^4)\otimes E[y_{11}, z_{13}]$.  In particular $P_*$ is a free differential graded $H^*(DI(4))$-module. Furthermore,  as a chain complex it is split as the tensor product of the following acyclic subcomplexes
\[\{P[x_7]/(x_7^4)\otimes P[\hat{a}_6]/(\hat{a}_6^4)\otimes E[z_{13}]\otimes\Gamma[\hat{t}_{24}, \hat{e}_{26}]\}\otimes \{E[y_{11}]\otimes \Gamma[\hat{b}_{10}]\}.\] Hence $P_*$ is a free $H^*(DI(4))$-resolution of $\F_2$.

Since $P_*$ is a free differential graded $H^*(DI(4))$-module, it is immediate that 
\begin{multline}\label{loop-space-DI(4)}E^2 = \Ext_{H^*(DI(4))}(\F_2,\F_2) = H^*(\Hom_{H^*(DI(4))}(P_*,\F_2))\cong \\ \Hom_{\F_2}(P[\hat{a}_6]/(\hat{a}_6^4)\otimes \Gamma[\hat{b}_{10}, \hat{t}_{24}, \hat{e}_{26}], \F_2)\cong P[a_6]/(a_6^2)\otimes P[b_{10}, c_{12}, e_{26}].\end{multline}
Since this module is concentrated in even degrees, there are no possible differentials,  so $E^2=E^\infty$. By inspection of the SSS for the path-loop fibration over $DI(4)$, one easily obtains $Sq_*^4(b_{10}) = a_6$ and $Sq_*^2(c_{12}) = b_{10}$. 

To calculate further Steenrod operations, we use a similar trick to the one used to $\Omega G_2$. Let $X$ denote the 7-connected cover of $DI(4)$. Thus there is a fibration 
\[X\rTo DI(4) \rTo{x_7} K(\Z, 7).\]
To calculate $H^*(X)$ we use Smith's Big Collapse Theorem \cite{Sm}. Notice that $x_7^*$ is onto, and its kernel consists of the ideal generated by all the polynomial generators of $H^*(K(\Z,7))$, different from $\iota_7$, $Sq^4\iota_7$, and $Sq^{6}\iota_7$, along with $\iota_7^4$, $(Sq^{6}\iota_7)^2$, and $(Sq^4\iota_7)^2$. This collection of generators forms a regular sequence in $H^*(K(\Z,7))$, and so the conditions of Smith's theorem are satisfied, and $H^*(X)$ can be written additively as the exterior algebra on infinitely many generators, corresponding in a 1-1 fashion to generators listed above, but with a dimension shift one down. Let $\epsilon_{27}$, $\tau_{25}$, and $\sigma_{21}$ be the elements in $H^*(X)$ corresponding to  $\iota_7^4$, $(Sq^{6}\iota_7)^2$, and $(Sq^4\iota_7)^2$ respectively.  Notice that $Sq^2(\tau_{25}) = \epsilon_{27}$ and $Sq^4(\sigma_{21}) = \tau_{25}$. Write
\begin{equation}\label{X}H^*(X) \cong E[\sigma_{21}, \tau_{25},\epsilon_{27}]\otimes E[k^2, k^3,\ldots, k^I,\ldots],\end{equation}
where $k^I$ in dimension $|Sq^I\iota_7|-1$ stands for the exterior generator corresponding to $Sq^I\iota_7$, for each $I$ such that $Sq^I\iota_7\in\Ker(x_7^*)$.

Now, consider the cohomology SSS for the principal fibration
\[\Omega DI(4)\rTo{j} K(\Z,6)\rTo{\pi} X.\]
Notice first that by naturality of the spectral sequence, the second factor in (\ref{X}) injects into $H^*(K(\Z,6))$ via $\pi^*$, while the classes $\sigma_{21}$, $\tau_{25}$ and $\epsilon_{27}$ are all in $\Ker(\pi^*)$. Furthermore, one has $j^*(\iota_6) = a_6$, and so  $j^*(Sq^4\iota_6) = \hat{b}_{10}$. The bottom dimensional class in $H^*(\Omega DI(4))$ which is not hit by $j^*$, is $\gamma_2(\hat{b}_{10})$, which is therefore transgressive. Hence $d(\gamma_2(\hat{b}_{10}))=\sigma_{21}$, and it follows that $d(Sq^4\gamma_2(\hat{b}_{10}))= Sq^4(\sigma_{21})=\tau_{25} $, while $d(Sq^6\gamma_2(\hat{b}_{10})) = d(Sq^{2,4}\gamma_2(\hat{b}_{10}))=Sq^2\tau_{25} = \epsilon_{27}$. 
 Hence $Sq^4\gamma_2(\hat{b}_{10}) = \hat{t}_{24}$, and $Sq^2\hat{t}_{24} = \hat{e}_{26}$. Dually, in homology, $Sq^2_*(e_{26}) = c_{12}^2$, and $Sq^4_*(c_{12}^2) = b_{10}^2$.

Next, we work out the Pontryagin algebra structure. Since $b_{10}$, $c_{12}$ and $e_{26}$ are elements of infinite height in $E^\infty$, they represent elements of infinite height in homology. Hence, it remains only to check whether $a_6^2=c_{12}$. But in cohomology one has $\hat{c}_{12} = Sq^2(\hat{b}_{10}) = Sq^2Sq^4(\hat{a}_6) = \hat{a}_6^2$. Hence in homology $\widebar{\Delta}(c_{12}) = a_6\otimes a_6$. But since $a_6$ is primitive, so is $a_6^2$, and since $H_{12}(\Omega DI(4))$ is 1-dimensional, it follows that $a_6^2 = 0$. This completes the calculation of the algebra structure.

The classes $a_6$ and $b_{10}$ are primitive for dimension reasons, and we have already computed the reduced diagonal of $c_{12}$.  Thus it remains to compute the reduced diagonal of $e_{26}$.  Notice that $H_{26}(\Omega DI(4))$ is 2-dimensional, generated additively by  $e_{26}$ and $a_6b_{10}^2$, and that $e_{26}$ can be modified by an additive summand of $a_6b_{10}^2$ without changing the algebra structure. For any choice of $e_{26}$ one has 
\[\widebar{\Delta}(e_{26}) = A(a_6b_{10}\otimes b_{10}+b_{10}\otimes a_6b_{10}) + B(a_6\otimes b_{10}^2 + b_{10}^2\otimes a_6),\]
for some $A,B\in\F_2$.  But $B$ is the coefficient of $\widebar{\Delta}(a_6b_{10}^2)$, and so by modifying the choice of $e_{26}$ if necessary, we may assume that $B=0$. Furthermore, if $A\neq 0$, then $H_{26}(\Omega DI(4))$ contains no primitive class, and so dually every class in $H^{26}(\Omega DI(4))$ is decomposable, which is clearly impossible. Hence $A=0$, and there is a choice for the class $e_{26}$ which is primitive. 

This completes the calculation of $H_*(\Omega DI(4))$ as a Hopf algebra and hence the proof of the proposition.
\end{proof}

Dually, the cohomology Hopf algebra is given by 
\[H^*(\Omega DI(4)) = P[\hat{a}_6]/(\hat{a}_6^4) \otimes \Gamma[\hat{b}_{10}, \hat{t}_{24}, \hat{e}_{26}].\]

We are now ready to start the calculation of $H_*(\Omega B\Sol(q))$. Consider the fibration
\[\Omega DI(4)\rTo \Omega B\Sol(q)\rTo DI(4).\]
The homology SSS associated to this fibration is a spectral sequence of Hopf algebras over $H_*(\Omega DI(4))$, whose $E^2$-page has the form
\begin{multline*}
E^2_{*,*} = H_*(\Omega DI(4))\otimes H_*(DI(4)) \cong\\ 
\left\{P[a_6]/(a_6^2) \otimes P[b_{10}, c_{12}, e_{26}]\right\}\otimes \left\{E[y_7, y_{11}, y_{13}]\otimes P[y_{14}]/(y_{14}^2)\right\}.\end{multline*}
Since $B\Sol(q)$ is 6-connected, $\pi_6(\Omega B\Sol(q)) \cong H_7(B\Sol(q),\Z)$. Thus $H_7(B\Sol(q),\Z)$ is a finite 2-group,  and so $d_7(y_7)=0$. Thus $d_7$ vanishes on   all the $y_i$ ($y_{14}$ by considering the dual cohomology spectral sequence, and the other generators by dimension reasons). The next possible nonvanishing differential is $d_{11}$. Considering the dual cohomology SSS, $d_{11}(\hat{b}_{10}) = d_{11}(Sq^4\hat{a}_6) = Sq^4d_7(\hat{a}_6) = 0$. Hence in homology $d_{11}(y_{11}) = 0$. Similarly $d_{13}(y_{13}) = 0$. Hence the spectral sequence collapses at $E^2$, and 
\begin{equation}\label{Solq-alg1}
H_*(\Omega B\Sol(q)) \cong P[a_6]/(a_6^2) \otimes P[b_{10}, c_{12}, e_{26}]\otimes E[y_7, y_{11}, y_{13}]\otimes P[y_{14}]/(y_{14}^2)\end{equation}
as a module over $H_*(\Omega DI(4))$.

\newcommand{\rd}{\widebar{\Delta}}
The loop space homology $H_*(\Omega DI(4))$ is contained in $H_*(\Omega B\Sol(q))$ as a Hopf subalgebra. The classes $y_7$ and $y_{11}$ are primitive for dimension reason. For $y_{13}$ one has $\rd(y_{13}) = A(a_6\otimes y_7) + y_7\otimes a_6) = A\rd(a_6y_7)$ for some $A\in\F_2$. Hence $y_{13}$ can be chosen to be primitive. Finally $\rd(y_{14}) = y_7\otimes y_7$. This completes the description of the coalgebra structure on $H_*(\Omega B\Sol(q))$.

Since $y_7$ and $y_{11}$ and $y_{13}$ are primitive, so are their squares. Since there are no nontrivial primitives in the respective dimensions, except for $e_{26}$, we conclude at once that $y_7 ^2 = y_{11}^2=0$, while $y_{13}^2 = Ae_{26}$ for some $A\in \F_2$.  But $Sq_*^2(e_{26})  = c_{12}^2$, while $Sq_*^2(y_{13}^2) = 0$, so $A=0$, and therefore $y_{13}^2 = 0$. Finally, $y_{14}^2$ is  primitive, and since there are no nontrivial primitives in dimension 28, it follows that $y_{14}^2=0$. 

Next, we calculate all the commutators involving the classes $y_i$. The results are summarised in the following table, while the calculations are below. Each entry in the table stands for the commutator $[\mathrm{Column}, \mathrm{Row}]$.

\begin{center} 
\begin{tabular}{||l|c|c|c|c|c|c|c|l||} \hline
               &$a_6$&$y_7$&$b_{10}$&$y_{11}$&$c_{12}$&$y_{13}$&$y_{14}$ & $e_{26}$ \\\hline
$y_7$  	&     $0$&$0$   &    $0$      &      $0$  &    $0$      &      $0$  &      $0$   & $0$\\\hline
$y_{11}$&     $0$&$0$   &       $0$   &      $0$  &    $0$     &       $0$  &      $0$   & $0$\\\hline
$y_{13}$&     $0$&$0$   &      $0$    &      $0$  &    $0$     &       $0$  &      $0$   & $0$\\\hline
$y_{14}$&     $b_{10}^2$&$0$   &        $c_{12}^2$ &       $0$  &       $e_{26}+a_6b_{10}^2$ &       $0$  &      $0$   & $b_{10}^4$\\\hline
\end{tabular}
\end{center}

Since $a_6$ and $y_7$ commute, $[a_6,y_{14}]$ is primitive, and so must be a multiple of $b_{10}^2$. Applying successive dual Steenrod squares
\[[a_6,y_{14}]\rMapsto{Sq_*^1}\left[a_6,y_{13}\right]\rMapsto{Sq_*^2} \left[a_6,y_{11}\right]\rMapsto{Sq_*^4}\left[a_6,y_7\right],\]
we conclude that all these commutators, with the possible exception of $[a_6, y_{14}]$ itself, vanish.

The class $y_7$ clearly commutes with itself, and its commutators with all other $y_i$ are primitive. This implies at once that $[y_7,y_{11}]$ and $[y_7, y_{14}]$ vanish, and $[y_7, y_{13}] = Sq^1_*[y_7,y_{14}] = 0$ as well.

The class $b_{10}$ commutes with $y_7$ for dimension reasons, and so $[b_{10}, y_{14}]$ is primitive, and is therefore a multiple of $c_{12}^2$. Applying dual Steenrod squares we have 
\[[b_{10},y_{14}]\rMapsto{Sq_*^1}\left[b_{10},y_{13}\right]\rMapsto{Sq_*^2} \left[b_{10},y_{11}\right],\]
which show that $[b_{10},y_{13}]$ and $[b_{10},y_{11}]$ vanish.

Since $y_{11}$ commutes with $y_{7}$, the commutator $[y_{11}, y_{14}]$ is primitive and $Sq_*^1[y_{11}, y_{14}] =[y_{11}, y_{13}]$. But there are no nontrivial primitives in dimension 25, and so both commutators vanish.

The classes $c_{12}$  and $y_7$ commute, since there are no 19 dimensional nonzero primitives, and so $[y_{14}, c_{12}]$ is primitive. Thus $[y_{14}, c_{12}]$ is a multiple of $e_{26}$. As before, we have
\[[c_{12},y_{14}]\rMapsto{Sq_*^1}\left[c_{12},y_{13}\right]\rMapsto{Sq_*^2} \left[c_{12},y_{11}\right],\]
which show that all commutators involving $c_{12}$, except possibly $[c_{12},y_{14}]$ vanish.

We already established that $y_{13}$ commutes with $y_7$ and $y_{11}$, and it commutes with $y_{14}$ as well since the commutator $[y_{13},y_{14}]$ is primitive. This also shows that all commutators with $y_{14}$ with other $y_i$ vanish. 

Finally, $[e_{26}, y_7]$ vanishes for lack of primitives in dimension 33. Thus $[e_{26}, y_{14}]$ is primitive and one has a chain of operations
\[\left[e_{26}, y_{14}\right]\rMapsto{Sq_*^1} \left[e_{26}, y_{13}\right] \rMapsto{Sq_*^2} \left[e_{26}, y_{11}\right].\]
The only nonzero primitive in dimension 40 is $b_{10}^4$, and so $[e_{26}, y_{14}]  =Ab_{10}^4$ for some $A\in \F_2$. Applying $Sq_*^1$ and $Sq_*^{2,1}$ to $b_{10}^4$, we conclude that $[e_{26}, y_{13}]$ and $[e_{26}, y_{11}]$ vanish.

It remains to evaluate the commutators of $y_{14}$ with the generators of $H_*(\Omega DI(4))$. To do that, we consider the cobar spectral sequence for $H_*(\Omega B\Sol(q))$, with
\[E^2 = \Cotor^{H_*(B\Sol(q))}(\F_2,\F_2) \cong H_*(T(\Sigma^{-1}\widebar{H}_*(B\Sol(q))),\, d_E),\]
where $d_E$ is the external differential on the cobar construction, induced by the reduced diagonal in $H_*(B\Sol(q))$.

Recall from \cite{Grbic}
\[H^*(B\Sol(q)) \cong
P[u_8, u_{12}, u_{14}, u_{15}, t_7, t_{11}, t_{13}]/I,\]
where $I$ is the ideal generated by the polynomials
$r_1=t_{11}^2 + u_8t_7^2+u_{15}t_7$, $r_2=t_{13}^2 + u_{12}t_7^2 + u_{15}t_{11}$ and $r_3=t_7^4 + u_{14}t_7^2 + u_{15}t_{13}$.

As for $BG_2(q)$, we denote classes in $H_*(B\Sol(q))$ by its dual cohomology class decorated by a bar. If $a\in \widebar{H}^*(B\Sol(q))$ is any class, then the corresponding tensor algebra generator will be denoted by $[\widebar{a}]$, while products of these generators will be written using the usual bar notation $[\widebar{a_1}|\widebar{a_2}|\cdots|\widebar{a_n}]$. Thus the relation $r_1$ translate to the following equation in the $E^1$ page of the cobar spectral sequence.
\[0 = d_E([\widebar{t_{11}^2}] + [\widebar{u_8t_7^2}] + [\widebar{u_{15}t_7}]) = [\widebar{t_{11}}]^2 + d_E( [\widebar{u_8t_7^2}]) + [[\widebar{u_{15}}],[\widebar{t_7}]].\]
The classes $[\widebar{t_{11}}]$, $[\widebar{u_{15}}]$ and $[\widebar{t_7}]$ are easily seen to be the permanent cycles in the spectral sequence corresponding to $b_{10}$, $y_{14}$ and $a_6$ respectively. Hence we obtain the relation
\[[a_6, y_{14}] = b_{10}^2.\]

Next, notice that one has 
\[\left[c_{12},y_{14}\right]\rMapsto{Sq_*^2}\left[b_{10}, y_{14}\right]\rMapsto{Sq_*^4}\left[a_6, y_{14}\right].\] The commutator $[b_{10},y_{14}]$ is primitive, while 
\[\widebar{\Delta}([c_{12},y_{14}]) = a_6\otimes [a_6,y_{14}] + [a_6,y_{14}]\otimes a_6 = a_6\otimes b_{10}^2 + b_{10}^2\otimes a_6.\]
 Hence we conclude that 
\[[b_{10}, y_{14}] = c_{12}^2 \quad\mathrm{and}\quad [c_{12},y_{14}] = e_{26} + a_6b_{10}^2 = e_{26} + a_6[a_6,y_{14}].\]
Finally, by the previous calculations,
\[[e_{26}, y_{14}] = [[c_{12},y_{14}], y_{14}]+ [a_6[a_6,y_{14}], y_{14}] = [a_6,y_{14}]^2 = b_{10}^4.\]
This completes the calculation of $H_*(\Omega B\Sol(q))$ as a Hopf algebra over the dual Steenrod algebra.

Our final task is the calculation of the BSS for $\Omega B\Sol(q)$. Using the known structure of $H^*(BDI(4))$ and $H^*(DI(4))$ and the corresponding BSS, we conclude that
\[H^*(BDI(4), \Z)\cong P[\upsilon_8, \upsilon_{12}, \upsilon_{15}, \upsilon_{28}]/(2\upsilon_{15}),\]
which allows us to conclude that $(\psi^q)^*(\upsilon_{2i}) = q^i\upsilon_{2i}$.
Also, similarly to the corresponding computation for $G_2$, 
\[H_i(DI(4), \Z) \cong  \begin{cases}
								\Z & i=0, 7, 11, 18, 27, 34, 38, 45\\
								\Z/2 & i = 13, 20, 24, 31
					\end{cases},\]
while all other homology groups vanish.  Denote homology classes by $\chi_i$, where $i$ corresponds to the dimension. The by inspection of the mod 2 homology structure, it is easy to conclude that $\chi_7$, $\chi_{11}$, and $\chi_{27}$ are the indecomposable among the torsion free classes, while the only torsion indecomposable class is $\chi_{13}$.

Consider the fibration
\[\Omega DI(4)\rTo \Omega B\Sol(q)\rTo DI(4).\]
The integral homology SSS calculation of this fibration is similar to the one done for $G_2$. 
Since $H_*(\Omega DI(4),\Z)$ is torsion free, the $E^2$ page of the spectral sequence is the tensor product of the homologies of base and fibre. 
One has a commutative diagram of fibrations
\[\begin{diagram}
\Omega DI(4) &\rTo& \Omega B\Sol(q) &\rTo& DI(4)\\
\dIgual && \dTo && \dTo_{f^q}\\
\Omega DI(4)& \rTo &* & \rTo&  DI(4)
\end{diagram}
\]
where $f^q$ is any map in the homotopy class of the composite
\[DI(4)\rTo{\Delta}DI(4)\times DI(4)\rTo{\Omega\psi^q\times (-1)} DI(4)\times DI(4) \rTo{\mu} DI(4).\] 
It is easy to verify that
\[H_n(DI(4)\times DI(4),\Z) \cong \bigoplus_{i+j=n}H_i(DI(4),\Z)\otimes H_j(DI(4),\Z)\]
for $n\le 32$. Hence it follows that for 
$i = 4, 6, 14$, one has $(f^q)_*(\chi_{2i-1}) = (q^i-1)\chi_{2i-1}$. 
In the SSS for the path-loop fibration over $DI(4)$, one has $d_7(\chi_7)= a_6$, $d_7(a_6\chi_7) = 2c_{12}$, $d_{11}(\chi_{11}) = b_{10}$, $d_{13}(\chi_{13}) = \bar{c}_{12}$ (the class of $c_{12}$ modulo $\Im(d_7)$), and $d_{27}(\chi_{27}) = e_{26}$.
Thus by commutativity of the diagram above, and naturality of the SSS one has in the integral homology SSS for the top row, $d_7(\chi_7) = (q^{4} - 1)$, $d_{11}(\chi_{11}) = (q^6-1)b_{10}$, and $d_{27}(\chi_{27}) =  (q^{14}-1)e_{26}$. Setting   $q = 4k\pm 1$ and $r_i = \nu_2(q^i-1)$, and performing the necessary arithmetics, we see that $r_6=r_{14} = \nu_2(k) +3 = r_4-1$.

We are now ready to compute the Bockstein spectral sequence for $H_*(\Omega B\Sol(q))$, which is a spectral sequence of modules over $H_*(\Omega DI(4))$, and  so we use the module structure given by Equations (\ref{Solq-alg1}) in the calculation. The first page of the spectral sequence is determined by $Sq^1_* (y_{14})=y_{13}$. Thus 
\[E^2 \cong P[a_6]/(a_6^2)\otimes P[b_{10}, c_{12}, e_{26}]\otimes E[y_7, y_{11},h_{27}],\]
where $h_{27}$ corresponds to the infinite cycle in $E^1$ given by $y_{13}y_{14}$.
By the integral homology calculation above, the next nontrivial differential is $\beta^{r_4-1}_*$, and one has $\beta^{r_4-1}(y_{11}) = c_{10}$, while $\beta^{r_4-1}(h_{27}) = e_{26}$. 
Next, we have $\beta_*^{r_4}(y_7) = a_6$, and since $\widebar{\Delta}(c_{12}) = a_6\otimes a_6$, it follows that $\beta^{r_4+1}_*(a_6y_7)$ is defined and is equal to $c_{12}$, provided that $a_6y_7\neq 0$ in $E^{r_4+1}$, which is obvious since already $E^2$ does not contain a nonzero class in dimension 14.  
This completes the calculation of the  BSS for $\Omega B\Sol(q)$, which  takes the form
\begin{center}
 \begin{tabular}{||l|c|c|c|c|c|c|c|c|l||} \hline
& $a_6$ & $y_7$ & $b_{10}$ & $y_{11}$  & $c_{12}$ &$y_{13}$& $y_{14}$ & $a_6y_7$ & $y_{13}y_{14}$ \\\hline
$Sq_*^1$ &$0$&$0$&$0$&$0$&$0$&$0$&$y_{13}$&$0$&$0$ \\ \hline
$\beta_*^{r_4-1}$ &$0$&$0$&$0$&$b_{10}$&$0$&$-$&$-$&$0$&$e_{26}$\\ \hline
$\beta_*^{r_4}$ &$0$&$a_6$&$-$&$-$&$0$&$-$&$-$&$0$&$-$\\ \hline
$\beta_*^{r_4+1}$ &$-$&$-$&$-$&$-$&$0$&$-$&$-$&$c_{12}$&$-$  \\ \hline
\end{tabular}
\end{center}

\end{document}